\documentclass[12pt,a4paper]{article}

\usepackage[utf8]{inputenc}
\usepackage[T1]{fontenc}
\usepackage{geometry}
\usepackage{amssymb}
\usepackage{amsmath}
\usepackage{amsthm}
\usepackage[all]{xy}
\usepackage{hyperref}
\usepackage{color}
\usepackage{mathrsfs}
\theoremstyle{plain}
\newtheorem{theorem}{Theorem}[section]
\newtheorem{proposition}[theorem]{Proposition}
\newtheorem{lemma}[theorem]{Lemma}
\newtheorem{corollary}[theorem]{Corollary}
\theoremstyle{definition}
\newtheorem{definition}[theorem]{Definition}

\newtheorem{remark}[theorem]{Remark}

\newcommand{\ra}{\rightarrow}
\newcommand{\lra}{\longrightarrow}

\newcommand{\Filt}[1]{\mbox{\rm Filt}(#1)}
\newcommand{\filt}[1]{\mbox{\rm filt}(#1)}

\newcommand{\X}{\mathcal{X}}

\newcommand{\Hom}{\operatorname{Hom}}
\newcommand{\Ext}{\operatorname{Ext}}
\newcommand{\End}{\operatorname{End}}
\renewcommand{\ker}{\operatorname{ker}}

\newcommand{\coker}{\operatorname{coker}}
\renewcommand{\mod}[1]{\operatorname{#1--mod}}
\newcommand{\Mod}[1]{\operatorname{#1--Mod}}

\newcommand{\pd}[1]{\mbox{\rm{proj.dim\,}}#1}

\title{Pr\"ufer modules in filtration categories of semibricks}

\author{Frank Lukas\\
\small Marc-Chagall-Str. 10, 51375 Leverkusen, Germany\\
\small Email: drfranklukas@gmail.com}

\date{}

\begin{document}

\maketitle

\begin{abstract}
	Let $R$ be a ring with unity. A brick in the module category $\Mod{R}$ is a finitely presented module whose endomorphism ring is a division ring. Two non-isomorphic bricks $X,Y$ are said to be orthogonal if $\Hom(X,Y)=\Hom(Y,X)=0$. A class $\X$ of pairwise orthogonal bricks is called a semibrick. We study the full subcategory $\Filt{\X}$ consisting of all modules admitting a filtration in $\X$ and show that this category is a wide subcategory of $\Mod{R}$.
	\\%
	Let $\X^{\perp}$ be the class of modules $M$ with $\Ext^1(X, M)=0$ for all $X\in \X$. For every module $Y\in \Mod{R}$ there exists an $\X^{\perp}$-envelope $Y_{\X}(\infty)$, which can be constructed as a direct limit of iterated universal short exact sequences. 
	\\%
	Assume in addition that every $X\in \X$ has projective dimension at most one. If moreover $\Hom(X, R)=0$ for all $X\in \X$, then the $\X^{\perp}$-envelope $R_{\X}(\infty)$ is isomorphic to the universal localization $R_{\X}$ of the ring $R$ at $\X$ in the sense of Schofield \cite{schofield1985representation}. We call the $\X^{\perp}$-envelopes of modules in $\X$ Pr\"ufer modules since they share many properties with Pr\"ufer groups and also with Pr\"ufer modules over tame hereditary algebras. Every injective object in $\Filt{\X}$ is isomorphic to a direct sum of Pr\"ufer modules.
\end{abstract}
\noindent\textbf{Mathematics Subject Classification (2020):} 16G60

\noindent\textbf{Keywords:} Pr\"ufer modules, representation theory of algebras, 
infinite dimensional modules, filtration categories, semibricks

\section*{\em In memory of my teacher Otto Kerner}
\section{Motivation and introduction}
Our motivation for analyzing filtration categories of semibricks was to construct Pr\"ufer modules over wild hereditary algebras. The starting point was the following observation.
\\%
The Pr\"ufer groups are precisely the injective objects in the abelian category of torsion groups. This category can be viewed as the filtration category of the semibrick $\X= \{\,\mathbb{Z}/p\mathbb{Z}\mid p \text{ prime}\,\}$. An analogous situation arises in the case of a finite-dimensional tame hereditary algebra. The module category of such an algebra possesses a maximal semibrick $\X$ that contains a representative of each isomorphism class of quasi-simple modules. The Pr\"ufer modules, as defined by Ringel in \cite[p. 326]{ringel1979infinite}, are the injective hulls within the wide subcategory of $\Mod{R}$ consisting of all regular modules. This category coincides with $\Filt{\X}$, the category of modules admitting an $\X$-filtration. 
\\%
The module category of a finite-dimensional wild hereditary algebra also has semibricks.  Let $X$ be a brick with self-extensions. By \cite[Theorem 3.14]{asai2025extensions}, there exists an infinite semibrick $\X$ containing $X$. Concrete examples arise from elementary modules, that is, nonzero regular modules that cannot occur as the middle term of a nontrivial short exact sequence of regular modules. As observed in \cite[Corollary 1.4(b)]{kerner1996elementary}, any class $\X$ of elementary modules sharing the same dimension vector forms a semibrick. In contrast to the tame case, there exist semibricks $\X$ with $dim_{\End(X)} \Ext^1(X, X) \geq 2$ and $\Ext^1(X, Y)\ne 0$ for all $X,Y\in \X$. 
\\%
The category $\filt{\X}$ of modules with an $\X$-filtration of finite length is a wide subcategory of $\mod{R}$, as Ringel proved in \cite[Theorem 1.2]{ringel1976representations}. In Section \ref{The_filtration_category_of_a_semibrick_is_abelian} of this paper, we generalize this result: the subcategory $\Filt{\X}$, consisting of all modules in $\Mod{R}$ that admit an $\X$-filtration of arbitrary length, is a wide subcategory of $\Mod{R}$. 
\\%
This result was first established by Claudius Heyer, who kindly shared with me his unpublished categorical proof.  I subsequently developed an algebraic proof, which is presented here. In the proof we exploit the fact that every module in $\Filt{\X}$ admits a canonical filtration, called the $\X$-socle filtration.  Moreover, every homomorphism between two modules in $\Filt{\X}$ is compatible with their $\X$-socle filtrations.  
\\%
In Section \ref{Universal_short_exact_sequences} we consider the $\X$-universal short exact sequence 
$$ \xymatrix{
	0 \ar[r] & M \ar[r] &  M_{\X}(2) \ar[r] & \bigoplus_{X\in \X} X^{(I_X) } \ar[r] & 0
} $$
which has the property that every short exact sequence in $\Ext^1(X,M)$ can be obtained by pull-backs for all $X\in \X$. This construction generalizes the approach introduced by Bongartz in \cite[Lemma 2.1]{Bongartz}. Iterating this process, we obtain a module $M_{\X}(\infty)$ as the direct limit of the ascending chain $M \subset M_{\X}(2) \subset M_{\X}(3) \subset \hdots$ of $\X$-universal inclusions. It turns out that the filtration $(M_{\X}(i), i\in \mathbb{N})$ is an $\X$-socle filtration whenever $M$ is a semisimple module in $\Filt{\X}$. From Lemma \ref{push-out_along_universal_inclusion} on, we assume that each brick in the semibrick $\X$ has projective dimension at most one. Under this  additional assumption we show that, for any module $M\in \Mod{R}$, the embedding $\alpha : M \ra M_{\X}(\infty)$ is an $\X^{\perp}$-envelope of $M$. 
\\%
In Section \ref{Pruefer_modules} we define a Pr\"ufer module as the $\X^{\perp}$-envelope of a brick in $\X$.  Pr\"ufer modules are indecomposable, and every module $M\in \Filt{\X}$ admits an injective hull $f : M\ra \hat{M}$ in $\Filt{\X}$,  where $\hat{M}$ is isomorphic to a direct sum of Pr\"ufer modules. 
\\%
I would like to thank Claudius Heyer for the fruitful discussions on the proof that $\Filt{\X}$ is a wide subcategory of $\Mod{R}$. Lidia Angeleri H\"ugel and the anonymous referee provided many valuable comments and corrections. They showed great patience with me as I returned to mathematics after a long period of working in different fields.
\section{Preliminaries}
Let $R$ be a ring with unity, and let $\Mod{R}$ (respectively $\mod{R}$) denote the category of all (respectively finitely presented) left $R$-modules. For a class $\X$ of modules in $\Mod{R}$ we denote by $Add(\X)$ the class of all modules  isomorphic to a direct summand of direct sums of modules of $\X$.
We write compositions from left to right: for $f\in \Hom(X,Y)$ and $g\in \Hom(Y,Z)$, the composition is $f\circ g\in \Hom(X,Z)$. A finitely presented module $X$ whose
 endomorphism ring is a division ring is called a {\bf brick}. Two non-isomorphic bricks $X$ and $Y$ are said to be orthogonal if $\Hom(X,Y)=\Hom(Y,X)=0$. A class of pairwise orthogonal bricks is called a {\bf semibrick}. In this paper we study filtration classes of semibricks. Part (a) of the following definition is taken from \cite[Definition 6.1]{gobel2012approximations}.
\begin{definition}
	
	Given a module $M$ and an ordinal  $\mu$, we call an ascending chain $\mathcal{M} = ( M_{\lambda}, \lambda \leq\mu )$ of submodules of $M$ a {\bf filtration of $M$}, if $M_0 = 0$, $M_{\mu} = M$, and $M_{\nu} = \bigcup_{\lambda < \nu} M_{\lambda}$ for every limit ordinal $\nu \leq \mu $. Moreover, for a class $\mathcal{C}$ of modules, we call $\mathcal{M}$ a {\bf $\mathcal{C}$-filtration of $M$} if each consecutive factor $M_{\lambda +1}/ M_{\lambda}$ ($\lambda \leq \mu$) is isomorphic to a module from $\mathcal{C}$. A module $M$ admitting a $\mathcal{C}$-filtration is said to be $\mathcal{C}$-filtered, and the class of all (finitely) $\mathcal{C}$-filtered modules is denoted by $Filt(\mathcal{C})$ (respectively $\filt{\mathcal{C}}$).
	\\%
	Let $M,N\in \Mod{R}$ have filtrations $( M_{\lambda}, \lambda \leq\mu )$ and $( N_{\lambda}, \nu \leq \omega )$. A homomorphism $f\in \Hom(M,N)$ is a {\bf homomorphism of the filtered objects }$ M$ and $ N$ if $f(M_{\lambda})\subset N_{\lambda}$ for all $\lambda \leq \mu$, where $N_{\lambda}=N$ for $\lambda > \mu$.
\end{definition}
For a class $\X$ we write $\Hom(\X,M)=0$ to mean $\Hom(X,M)=0$ for all $X\in \X$, and we define  
$$
\begin{array}{l}
	\X^{\perp} = \ker \Ext^1(\X, -) = \{ M\in \Mod{R} \, | \, \Ext^1(X,M)=0\ for\ all\  X\in \X\}.
\end{array}
$$
Similarly, the class ${}^{\perp}\X= \ker \Ext^1(-,\X)$ is defined in the dual way. From \cite[Lemma 6.2]{gobel2012approximations} we recall the  Eklof Lemma:
\begin{lemma} \label{eklof-lemma}
	Let $N$ be a module and $M$ be a ${}^{\perp}N$-filtered module. Then $M\in {}^{\perp}N$.
\end{lemma}
From \cite[Lemma 6.6]{gobel2012approximations} we also adopt the following result: 
\begin{lemma} \label{direct_limes_of_Ext}
	Let $R$ be a ring and $M$ be a finitely presented module. Let $\{ N_{\alpha}, f_{\beta, \alpha} | \alpha \leq \beta \in I\}$ be a direct system of modules. Then we have 
	$$ \Ext^1(M, \varinjlim_{\alpha \in I} N_{\alpha}) \cong \varinjlim_{\alpha\in I} \Ext^1(M, N_{\alpha})$$
\end{lemma}
For a module $M\in \Mod{R}$ we denote by $\pd{M}$ its projective dimension. For a class of modules $\X$ we write $\pd{\X}\leq 1$ if  $\pd{X}\leq 1$ for all $X\in \X$.
The following result, due to Auslander \cite[Lemma 6.4]{gobel2012approximations}, is useful.
\begin{lemma}\label{proj-dimension}
	Let $n\in \mathbb{N}$, and let $M$ be a module. Assume $M$ has a filtration of modules $X\in \X$ with $\pd{\X}\leq n$. Then $\pd{M}\leq n$. 
\end{lemma}
Finally, for modules $U_i, V \in \Mod{R}$ (indexed by $i\in I$) and $W\in \mod{R}$,  we recall the standard isomorphisms
$$ \Ext^1 (\bigoplus_{i\in I} U_i, V)\cong \prod_{i\in I} \Ext^1(U_i,V)  , \Ext^1 (W, \bigoplus_{i\in I} U_i)\cong \bigoplus_{i\in I} \Ext^1(W, U_i)$$
\section{The filtration category of a semibrick is a wide subcategory of $\Mod{R}$} \label{The_filtration_category_of_a_semibrick_is_abelian}
A full subcategory $\mathcal{C}$ of an abelian category $\mathcal{A}$ is called a {\bf wide subcategory of $\mathcal{A}$} if it is abelian and closed under extensions. For a semibrick $\X$, Ringel showed in \cite[Theorem 1.2]{ringel1976representations} that $\filt{\X}$ is a wide subcategory of $\mod{R}$. In this section we prove that $\Filt{\X}$ is a wide subcategory of $\Mod{R}$. 
\begin{definition}
	Let $\X$ be a semibrick. 
	\begin{itemize}
		\item[(a)]  For $M$ in $\Filt{\X}$ the {\bf $\X$-trace of $M$} is 
		$$tr_{\X}\ M = \sum_{f\in \Hom(X,M), X\in \X} Im\ f$$
		\item[(b)] A module $M\in \Filt{\X}$ is called {\bf $\X$-semisimple} if $M$ is isomorphic to  a direct sum of objects from $\X$. 
	\end{itemize}
\end{definition}
Note that for any nonzero module $M$ in $\Filt{\X}$, the submodule $\X$-trace of $M$ is nonzero. Every nonzero submodule $U\subseteq M$ with $U$ in $\Filt{\X}$ intersects $tr_{\X}\ M$ nontrivially; in this sense $tr_{\X}\ M$ is $\X$-essential in $M$.   
\begin{proposition} \label{trace_is_semisimple}
	Let $\X$ be a semibrick and $M\in \Filt{\X}$. The following statements are equivalent:
	\begin{itemize}
		\item[(i)] $M$ is $\X$-semisimple
		\item[(ii)] $ tr_{\X}\ M = M$
	\end{itemize}
\end{proposition}
	\begin{proof}
	$(i)\Rightarrow (ii)$: is clear.
	\\%
	$(ii) \Rightarrow (i)$: Let $(M_{\lambda}, \lambda \leq \mu)$ be an $\X$-filtration of $M$, and $\Lambda$ be the set of ordinals $0< \lambda \leq \mu$ for which $M_{\lambda -1 } \hookrightarrow M_{\lambda }$ splits. We claim that $tr_{\X} \ M \cong \bigoplus_{\lambda\in \Lambda} M_{\lambda } / M_{\lambda -1}$.
	\\%
	For $\lambda \in \Lambda$ the projection $M_{\lambda} \ra M_{\lambda}/ M_{\lambda -1}$ splits; let $r_{\lambda}: M_{\lambda}/ M_{\lambda -1} \ra M_{\lambda}$ be a retraction and set $U_{\lambda}= Im\ r_{\lambda}$. We show by transfinite induction that $tr_{\X} \ M_{\lambda} = \bigoplus_{\alpha \in \Lambda, \alpha\leq \lambda} U_{\alpha}$.
	\\%
	For $\lambda=1$ we have $M_1=U_1= tr_{\X}\ M_1$. Let now $\lambda$ be an ordinal. If $M_{\lambda -1} \hookrightarrow M_{\lambda}$ splits, then $M_{\lambda} = M_{\lambda -1} \oplus U_{\lambda}$ with $U_{\lambda}$  isomorphic to a brick $X \in \X$. Any homomorphism $f: X\ra M_{\lambda -1} \oplus U_{\lambda}$ with $X$ in $\X$ can be written as $(f_1, f_2)$ with $f_1\in \Hom(X, M_{\lambda -1})$ and $f_2\in \Hom(X,U_{\lambda})$. The image of $f_1$ lies in $tr_{\X} M_{\lambda -1} = \oplus_{\alpha\in \Lambda,\alpha\leq \lambda -1} U_{\alpha}$, which shows that $tr_{\X} M_{\lambda}= \oplus_{\alpha\in \Lambda,\alpha\leq \lambda} U_{\alpha}$.
	\\%
	 If $M_{\lambda -1} \hookrightarrow M_{\lambda}$ does not split, then for the projection $\pi:M_{\lambda} \ra M_{\lambda} / M_{\lambda-1}$ and any $f:X\ra M_{\lambda}$ with $X\in \X$ the composition $f\circ \pi$ is not an isomorphism and hence zero, so $Im \ f\subset M_{\lambda-1}$.
	 \\%
	 For a limit ordinal $\lambda$, $tr_{\X}\  M_{\lambda} = \bigcup_{\alpha\in \Lambda , \alpha< \lambda} tr_{\X}\  M_{\alpha}= \bigoplus_{\alpha \in \Lambda, \alpha < \lambda} U_{\alpha}$. Thus $M$ is the direct sum of the bricks $U_{\lambda}$ and hence $\X$-semisimple. 
	\end{proof}
By this proposition, every module in $Add(M)$ is $\X$-semisimple if and only if $M$ is $\X$-semisimple. For $M$ in $\Filt{\X}$ we write {$soc_{\X}\ M$  instead of $tr_{\X}\ M$ and call it the {\bf $\X$-socle of $M$}.
	\begin{proposition} \label{existence_X_socle_filtration}
		Let $\X$ be a semibrick and $M\in \Filt{\X}$ with an $\X$-filtration $(M_{\lambda}, \lambda \leq \mu)$. There exists a filtration $(N_{\lambda}, \lambda \leq \mu)$ of $M$ and $\xi\leq \mu$ such that $N_{\xi} = soc_{\X} M$. In particular,   $M / soc_{\X} M \in \Filt{\X}$. 
	\end{proposition}
	\begin{proof} 
		With $\Lambda$ as in the previous proof, for each $\lambda \in \Lambda$ there is a submodule $U_{\lambda}$ of $M_{\lambda}$ such that $M_{\lambda} = U_{\lambda} \oplus M_{\lambda -1}$ and $soc_{\X} \ M = \bigoplus_{\lambda \in \Lambda} U_{\lambda}$. 
		\\%
		Let $(W_{\lambda}, \lambda \leq \xi)$ be an $\X$-filtration of $soc_{\X}\ M$. Define a filtration $(N_{\lambda}, \lambda \leq \mu + \xi )$ by $N_{\lambda}= W_{\lambda}$ for $\lambda \leq \xi$, $N_{\xi + \lambda} = soc_{\X} M + M_{\lambda}$ for $ 0 < \lambda \leq \mu$.
		We claim that for an ordinal  $\lambda$ either $N_{\lambda}$ is equal to $N_{\lambda -1}$ or $N_{\lambda}/ N_{\lambda -1}$ is isomorphic to a brick in $\X$. 
		If $\lambda$ is lower or equal $\xi$ there is nothing to prove. 
		Therefore we consider an ordinal of the form $\xi + \lambda$ with $0 < \lambda \leq \mu$. If $\lambda \in \Lambda$ then we have 
		$$ N_{\xi + \lambda} = soc_{\X} M + M_{\lambda}  = soc_{\X} M + ( U_{\lambda}  \oplus M_{\lambda -1})  = N_{(\xi + \lambda) -1}$$ 
		since $U_{\lambda} \subset soc_{\X} M$. It remains to consider $N_{\xi + \lambda} / N_{(\xi + \lambda) -1}$ for $0< \lambda \leq \mu$ with $\lambda \notin \Lambda$. 
		In this case we have:
		$$
		\begin{array}{rcl}
			N_{\xi + \lambda} / N_{(\xi + \lambda) -1} & = & ( soc_{\X} M + M_{\lambda}) / ( soc_{\X} M + M_{\lambda -1 }) \\
			& = & ( \bigoplus_{\omega \geq \lambda, \omega \in \Lambda} U_{\omega} \oplus M_{\lambda }) / ( \bigoplus_{\omega \geq \lambda, \omega \in \Lambda} U_{\omega} \oplus M_{\lambda -1} ) \\
			& \cong & M_{\lambda } / M_{\lambda -1} \in \X
		\end{array}
		$$
		For a limit ordinal $\xi + \lambda$ with $0< \lambda \leq \mu$ we have 
		$$N_{\xi + \lambda}  = soc_{\X} M + M_{\lambda}   = \bigcup_{\alpha < \lambda} (soc_{\X}  M +M_{\alpha} )   = \bigcup_{\alpha < \lambda} N_{\xi + \alpha} =\bigcup_{\alpha <  \xi + \lambda} N_{\alpha} $$	
		If we remove from the filtration $(N_{\lambda}, \lambda \leq \mu + \xi )$ the $N_{\lambda}$ with $N_{\lambda}= N_{\lambda +1}$, then we obtain an $\X$-filtration of $M$ with the properties $N_{\xi} = soc_{\X}\ M$ and $N_{\mu + \xi}=M$.
	\end{proof}
	Part (a) of the following lemma can be viewed as an analogue of Steinitz's theorem. In part (b) we use the elementary observation that any homomorphism $f \in \Hom(X,M)$ with $X \in \X$ and $M \in \Filt{\X}$ is either injective or zero.
	\begin{lemma}\label{intersections_of_modules}
		Let $\X$ be a semibrick and let $M$ be an $\X$-semisimple module. Suppose $M=\bigoplus_{i \in I} M_i$, where each $M_i$ is isomorphic to a brick in $\X$.		
		\begin{itemize}
			\item[(a)] Let $U$ be a finitely presented $\X$-semisimple submodule of $M$. Then there exists a subset $I' \subseteq I$ such that $M=\bigoplus_{i \in I'} M_i \oplus U$.			
			\item[(b)] Let $V$ be a submodule of $M$ which is isomorphic to a brick in $\X$, and let $U$ be an $\X$-semisimple submodule of $M$. If $V\cap U\neq 0$, then $V\subseteq U$.
		\end{itemize}
	\end{lemma}
	\begin{proof}
		(a) Since $U$ is finitely presented and $\X$-semisimple, we can write $U=\bigoplus_{j=1}^{r} U_j$ for some $r\in\mathbb N$, where each $U_j$ is isomorphic to a brick in $\X$. We prove the statement by induction on $r$. If $U=0$, there is nothing to prove. Assume the statement holds for $r\ge 1$ and let $U=\bigoplus_{j=1}^{r+1} U_j$. Set $U'=\bigoplus_{j=1}^{r} U_j$. By the induction hypothesis there exists $I'\subseteq I$ such that $M=\bigoplus_{i\in I'} M_i \oplus U'$. Since $U_{r+1}\cap U'=0$, there exists $i'\in I'$ such that the canonical projection  $\pi_{i'}:\bigoplus_{i\in I'} M_i \to M_{i'}$ satisfies $\pi_{i'}(U_{r+1})\neq 0$. As $U_{r+1}$ is isomorphic to a brick in $\X$, the restriction $\pi_{i'}|_{U_{r+1}}:U_{r+1}\to M_{i'}$ is an isomorphism. Hence $\bigoplus_{i\in I'} M_i=\ker \pi_{i'} \oplus U_{r+1}$ and therefore $M=\bigoplus_{i\in I'\setminus\{i'\}} M_i \oplus U$.
		\\%
		(b) Since $U$ is $\X$-semisimple, we can write $U=\bigoplus_{j\in J} U_j$. Let $0\ne v\in V\cap U$. Then $v\in\bigoplus_{j\in J'} U_j$ for some finite $J'\subseteq J$. Set $U'=\bigoplus_{j\in J'} U_j$. By part (a) there exists $I'\subseteq I$ such that $M=\bigoplus_{i\in I'} M_i \oplus U'$. Let $\pi:\bigoplus_{i\in I'} M_i \oplus U' \to \bigoplus_{i\in I'} M_i$ be the canonical projection. Since $V$ is isomorphic to a brick in $\X$, the map $\pi|_V$ is either zero or injective. Since $\pi(v)=0$, it follows that $\pi|_V=0$, and hence $V\subseteq\ker\pi=U'\subseteq U$.
	\end{proof}
	\begin{proposition} \label{semisimple_modules}
		Let $\X$ be a semibrick, let $N$ be a module with an $\X$-filtration, and let $M$ be an $\X$-semisimple module. For any homomorphism $f\in \Hom(N,M)$, the image of $f$ is an $\X$-semisimple submodule and a direct summand of $M$.
	\end{proposition}
	
	\begin{proof}
		Let $(N_{\lambda})_{\lambda\le \mu}$ be an $\X$-filtration of $N$. We proceed by transfinite induction on $\mu$.		
		If $\mu=0$, the statement is trivial. Assume that $\mu$ is an ordinal such that $f(N_{\mu-1}) \subsetneqq f(N_{\mu})$. By the induction hypothesis we can write $M=C\oplus f(N_{\mu-1})$. Since $C$ is $\X$-semisimple, we have $C=\bigoplus_{i\in I} C_i$ with each $C_i$ isomorphic to a brick in $\X$. Let $\pi_C : C\oplus f(N_{\mu-1}) \to C$ denote the canonical projection. Because $f(N_{\mu-1})$ is strictly contained in $f(N_{\mu})$, the module $\pi_C(f(N_{\mu}))$ is nonzero.
		
		We claim that $\pi_C(f(N_{\mu}))$ is isomorphic to a brick in $\X$. Consider the commutative diagram
		\[
		\xymatrix{
			N_{\mu} \ar[r]^f \ar[d] & f(N_{\mu}) \ar[r]^{\pi_C} \ar[d]^{\pi} & C \\
			N_{\mu}/N_{\mu-1} \ar[r]^{\bar f} & f(N_{\mu})/f(N_{\mu-1}) \ar[ur]^{\bar\pi_C}
		}
		\]
		in which the canonical projection $N_{\mu}\to N_{\mu}/N_{\mu-1}$ factors through $\bar f : N_{\mu}/N_{\mu-1} \to f(N_{\mu})/f(N_{\mu-1})$ and $\pi_C$ factors through $\bar\pi_C$. Since the composition $\bar\pi_C\circ \bar f$ is nonzero and $N_{\mu}/N_{\mu-1}$ is a brick, it follows that this map is a monomorphism. Hence $\pi_C(f(N_{\mu}))$ is isomorphic to a brick in $\X$.
		
	By Lemma~\ref{intersections_of_modules}(a) there exists a $I'\subseteq I$ with 
		$$	C=\bigoplus_{i\in I'} C_i \oplus \pi_C(f(N_{\mu})).$$
		Consequently
	$$ 	M=\bigoplus_{i\in I'} C_i \oplus \pi_C(f(N_{\mu})) \oplus f(N_{\mu-1}) =\bigoplus_{i\in I'} C_i \oplus f(N_{\mu}),$$
		so $f(N_{\mu})$ is a direct summand of $M$.
		
		Now assume that $\mu$ is a limit ordinal. By induction, each $f(N_{\lambda})$ for $\lambda<\mu$ is $\X$-semisimple. Since
		$$tr_{\X}f(N_{\mu}) =  \bigcup_ {\lambda < \mu}  tr_{\X} f(N_{\lambda}) = \bigcup_{\lambda < \mu} f(N_{\lambda})  = f(N_{\mu})$$
		Proposition~\ref{trace_is_semisimple} implies that $f(N_{\mu})$ is $\X$-semisimple.
		
		If $f$ is an epimorphism, the statement is trivial. Hence assume that $f(N_{\mu})$ is a proper submodule of $M$. Consider all families $(V_{\omega})_{\omega\in\Omega}$ of submodules of $M$ that are isomorphic to bricks in $\X$ such that $\sum_{\omega\in\Omega} V_{\omega} + f(N_{\mu})$ is direct. Write $M=\bigoplus_{i\in I} M_i$ with each $M_i$ isomorphic to a brick in $\X$. Since $f(N_{\mu})$ is a proper submodule of $M$, there exists $i\in I$ such that $M_i\nsubseteq f(N_{\mu})$. By Lemma~\ref{intersections_of_modules}(b) we have $M_i\cap f(N_{\mu})=0$. Thus the collection of such families is nonempty.
		
		For two families $(V_{\omega})_{\omega\in\Omega}$ and $(W_{\omega})_{\omega\in\Omega'}$ we define $(V_{\omega})_{\omega\in\Omega}\le (W_{\omega})_{\omega\in\Omega'}$ if $\Omega\subseteq\Omega'$ and $V_{\omega}=W_{\omega}$ for all $\omega\in\Omega$. This defines a partial order on the collection of these families. Let $\{(V_{\omega}^{\lambda})_{\omega\in\Omega_{\lambda}}\}_{\lambda\in\Lambda}$ be an ascending chain. Then the union of the index sets again yields a family which serves as an upper bound. Hence, by Zorn's lemma, there exists a maximal family $(V_{\omega})_{\omega\in\Omega_m}$.
		
		If $\bigoplus_{\omega\in\Omega_m} V_{\omega} \oplus f(N_{\mu}) \neq M$, then there exists $i\in I$ such that $M_i$ is not contained in this sum. By Lemma~\ref{intersections_of_modules}(b) we obtain
		$M_i \cap \left(\bigoplus_{\omega\in\Omega_m} V_{\omega} \oplus f(N_{\mu})\right)=0$,
		contradicting the maximality of $(V_{\omega})_{\omega\in\Omega_m}$. Hence
		$$	M=\bigoplus_{\omega\in\Omega_m} V_{\omega} \oplus f(N_{\mu}),	$$ 
		so $f(N_{\mu})$ is a direct summand of $M$.
	\end{proof}
\begin{proposition} \label{semisimple}
	Let $\X$ be a semibrick and $f: M \ra N$ be a nonzero homomorphism between two $\X$-semisimple modules $M$ and $N$. Then $\ker \ f$ and $\coker \ f $ are $\X$-semisimple.
\end{proposition}
\begin{proof}
	By Proposition \ref{semisimple_modules}, we have $N\cong Im \ f \oplus N'$ with an $\X$-semisimple module $N'$. So $\coker  f \cong N'$ is $\X$-semisimple. 
	\\%
	To show that $\ker  f$ is $\X$-semisimple, assume that $\ker  f \ne 0$. It suffices to show that the short exact sequence 
	\begin{equation} \label{splitting-sequence}
		\xymatrix{
			0 \ar[r] & \ker f \ar[r]^{\varepsilon} & M \ar[r]^{f} & Im f \ar[r] & 0
		}
	\end{equation}
	splits. 
	Note that $Im \ f$ is $\X$-semisimple by Proposition \ref{semisimple_modules}. Write $M=\bigoplus_{i\in I} M_i$, where each $M_i$ is isomorphic to a brick in $\X$. Consider the partially ordered set $({\mathcal{P}}, \leq )$ whose elements are subsets $I'\subseteq I$ such that the restriction of $f$ to  $\bigoplus_{i\in I'} M_i$ is injective. The set $\mathcal{P}$ is partially ordered by inclusion and non-empty, since there exists $i\in I$ such with $f(M_i)\ne 0$. Every chain in $\mathcal{P}$ has an upper bound (given by the union of the subsets), so by Zorn's Lemma there exists a maximal subset $I_m$ of $I$ such that $f$ restricted to $\bigoplus_{i\in I_m} M_i$ is injective. If $f(\bigoplus_{i\in I_m} M_i) = Im\ f$, then the short exact sequence \eqref{splitting-sequence} splits and we are done. Otherwise, by Proposition \ref{semisimple_modules} there exists a submodule $V\subseteq Im\ f$ such that $Im\ f = f(\bigoplus_{i\in I_m} M_i) \oplus V$. Let $\pi: Im\ f \ra V$ be the canonical projection. There exists a $i_0\in I \setminus I_m$ such that $f\circ \pi$ restricted to $M_{i_0}$ is injective - a contradiction to the maximality of $I_m$. Hence the short exact sequence \eqref{splitting-sequence} splits, and  $\ker  f$ is $\X$-semisimple. 
\end{proof}
Given a semibrick $\X$, an $\X$-filtration of a module $M\in \Filt{\X}$ is in general not uniquely determined. Using Proposition \ref{existence_X_socle_filtration}, we can define  a canonical filtration for every module in $\Filt{\X}$.
\begin{definition}
	Let $\X$ be a semibrick and $M$ be an $\X$-filtered module. The {\bf $\X$-socle filtration} $(V_\lambda, \lambda \leq \mu)$ of $M$ is defined inductively as follows: 
	\\%
	Start with $V_0=0$. If $V_\lambda$ is defined, let $V_{\lambda +1}$ be the pull-back along the inclusion $\varepsilon: soc_{\X} \ (M/V_{\lambda}) \hookrightarrow M/V_{\lambda}$ in the short exact sequence 
	$$
	\xymatrix{
		0 \ar[r] & V_{\lambda} \ar[r] & M \ar[r] & M / V_{\lambda} \ar[r] & 0 \\
		0 \ar[r] & V_{\lambda} \ar[r] \ar@{=}[u]  & V_{\lambda +1} \ar[r] \ar[u] & soc_{\X} (M / V_{\lambda}) \ar[r] \ar[u]_{\varepsilon} & 0
	}
	$$
		If $\lambda$ is a limit ordinal,  define $V_{\lambda} = \bigcup_{\alpha < \lambda} V_{\alpha}$.
	\end{definition}
	An $\X$-socle filtration $(M_{\lambda}, \lambda \leq \nu)$ of $M$ is a filtration with successive factors $M_{\lambda}/ M_{\lambda-1} \in Add(\X)$. Every module admitting an $Add(\X)$-filtration belongs to $\Filt{\X}$. 
	\begin{proposition} \label{strict-filtration}
		Let $\X$ be a semibrick, and let $M,N\in \Filt{\X}$. Then every homomorphism $f\in \Hom(M,N)$ is a homomorphism of the $\X$-socle filtered objects.
	\end{proposition}
	\begin{proof}
		Let $(M_{\lambda} ,  \lambda \leq \mu)$ and $(N_{\omega}, \omega \leq \nu)$ be the $\X$-socle filtration of $M$ and $N$, respectively. We must show that $f(M_{\lambda}) \subset N_{\lambda}$ for all $\lambda \leq \mu$. If $\lambda > \nu$, we define $N_{\lambda}=N$. 
		The proof proceeds by transfinite induction on $\lambda$. For $\lambda =0$, the claim is trivial. 
		Assume now that $f(M_{\lambda}) \subseteq N_{\lambda}$ for some $\lambda$, and let $f_{\lambda}$ denote the restriction of $f$ to $M_{\lambda}$. Consider the following commutative diagram, where $\varepsilon$ is the inclusion of $soc_{\X} \ (M/M_{\lambda}) $ into $M/M_{\lambda}$:
		$$ \xymatrix{
			0 \ar[r] & M_{\lambda} \ar[r] \ar@{=}[d]& M_{\lambda +1} \ar[r] \ar[d]& soc_{\X} (M / M_{\lambda} )\ar[r] \ar[d]^{\varepsilon} & 0 \\
			0 \ar[r] & M_{\lambda} \ar[r] \ar[d]^{f_{\lambda}} & M \ar[r] \ar[d]^{f}& M / M_{\lambda} \ar[r] \ar[d]^{\bar{f}} & 0 \\
			0 \ar[r] & N_{\lambda} \ar[r] & N \ar[r] & N/N_{\lambda} \ar[r] & 0
		} $$
		The image of $soc_{\X} M/M_{\lambda}$ under $\bar{f}$ is contained in $soc_{\X} (N/N_{\lambda})$. Hence 
		$f(M_{\lambda +1})  \subseteq N_{\lambda +1}$.  
		\\%
		If $\lambda$ is a limit ordinal, then $f(M_{\alpha}) \subseteq N_{\alpha}$ for all $\alpha < \lambda$ implies $$f(M_{\lambda}) = \bigcup_{\alpha < \lambda} f(M_{\alpha}) \subset \bigcup_{\alpha < \lambda} N_{\alpha}= N_{\lambda}.$$
		This completes the proof.
	\end{proof}
\begin{proposition} \label{closed_under_kernel}
	Let $\X$ be a semibrick. The category $\Filt{\X}$ is closed under extensions, kernels and images.
\end{proposition}
\begin{proof}
The argument that $\Filt{\X}$ is exact follows the reasoning of \cite[Lemma 3.2]{hugel2019characterisation}. Let $M,N\in \Filt{\X}$ have $\X$-socle filtrations $(M_{\lambda} , \lambda \leq \mu)$ and $(N_{\omega}, \omega \leq \nu)$, and consider a short exact sequence 
	$$ \xymatrix{
	0 \ar[r] & M \ar[r] & Z \ar[r] & N \ar[r] & 0
	} $$
For each $N_{\omega}$ with $\omega \leq \nu$, choose a subobject $V_{\mu + \omega} \subseteq Z$ such that $V_{\mu + \omega}/M \cong N_{\omega}$. Define $(Z_{\Theta}, \Theta\leq \mu + \nu)$ by setting $Z_{\Theta}= M_{\Theta}$ for $\Theta \leq \mu$ and $Z_{\Theta}= V_{\Theta}$ for $\mu < \Theta \leq \mu + \nu$. Then $(Z_{\Theta}, \Theta\leq \mu + \nu)$ is an $Add(\X)$-filtration of $Z$, since 
$$Z_{\mu +1}/ Z_{\mu} = Z_{\mu +1} / M = V_1/M \cong N_1 \in Add(\X) .$$
Hence $\Filt{\X}$ is closed under extensions. 
\\%
Now let $f:M\ra N$ be a homomorphism in $\Filt{\X}$. For each $\lambda \leq \mu$, denote by  $f_{\lambda}$ the restriction of $f$ to $M_{\lambda}$. By Proposition \ref{strict-filtration}, we have $f_{\lambda}(M_{\lambda}) \subseteq N_{\lambda}$ for all $\lambda$. Applying the Snake Lemma to the induced diagrams give exact rows: 
\begin{equation} \label{snake-diagram}
			\xymatrix{
			0 \ar[r]		&        \ker f_{\lambda -1}      \ar[r]    \ar[d]           & \ker f_{\lambda } \ar[r]^{\alpha} \ar[d]        &  \ker \bar{f_{\lambda }} \ar[r]^{\beta}  \ar[d] & \\
			0 \ar[r] & M_{\lambda -1} \ar[r] \ar[d]^{f_{\lambda -1}}  & M_{\lambda } \ar[r] \ar[d]^{f_{\lambda }} & M_{\lambda }/M_{\lambda -1} \ar[r] \ar[d]^{\bar{f_{\lambda }}} & 0 \\
			0 \ar[r] & N_{\lambda -1}  \ar[r] \ar[d] & N_{\lambda } \ar[r] \ar[d] & N_{\lambda } /N_{\lambda -1}  \ar[r] \ar[d]& 0 \\
			\ar[r]^{\beta\hspace{0.5cm}}	&        \coker f_{\lambda -1}      \ar[r]^{\gamma}            & \coker f_{\lambda } \ar[r]        &  \coker \bar{f_{\lambda }} \ar[r] & 0
		}
\end{equation}
	We show by transfinite induction on $\lambda$ that all modules in the upper row lie in $\Filt{\X}$. For $\lambda=1$ this follows directly from Proposition \ref{semisimple}. Assume the statement holds for all ordinals less than $\lambda$. Since both $M_{\lambda }/M_{\lambda -1}$ and $N_{\lambda } /N_{\lambda -1}$ belong to $Add(\X)$, Proposition \ref{semisimple} implies that $\ker \bar{f_{\lambda }}$ and $\coker \bar{f_{\lambda }}$ are also in $Add(\X)$. By the induction hypothesis, $\ker f_{\lambda -1}$ and $coker f_{\lambda -1}$ are in $\Filt{\X}$. 
	\\%
	Consider $\beta$ as a homomorphism from the $\X$-semisimple $\ker\ \bar{f_{\lambda }}$ to the $\X$-semisimple $\X$-socle of $\coker\ f_{\lambda-1}$. By Proposition \ref{semisimple}, $\ker\ \beta = Im\ \alpha$ is $\X$-semisimple, so $\ker f_{\lambda }\in \Filt{\X}$. 
	Since the image of $\beta$ lies in $soc_{\X} (\coker f_{\lambda-1})$, the quotient $\coker f_{\lambda -1}/ Im \ \beta $ is in $\Filt{\X}$, hence $\coker f_{\lambda} \in \Filt{\X}$.
	\\%
		If $\lambda$ is a limit ordinal, then $\ker f_{\lambda} = \bigcup_{\omega < \lambda} \ker f_{\omega}$, which shows $\ker f \in \Filt{\X}$.
	\\%
	To verify closure under images, note that $(f(M_{\lambda}), \lambda \leq \mu)$ forms an $Add(\X)$-filtration of $Im \ f$. For $\lambda =1$, Proposition \ref{semisimple_modules} gives $f(M_1)\in Add(\X)$. For an ordinal $\lambda$, define $f_{\lambda}'\in \Hom(M_{\lambda}, f(M_{\lambda}))$ as the restricted map and apply the snake lemma: 
	\begin{equation} \label{snake-diagram_2}
		\xymatrix{
			0 \ar[r]		&        \ker f'_{\lambda-1}      \ar[r]    \ar[d]           & \ker f'_{\lambda } \ar[r]^{\alpha} \ar[d]        &  \ker \bar{f'_{\lambda }} \ar[r]  \ar[d] & 0 \\
			0 \ar[r] & M_{\lambda-1} \ar[r] \ar[d]^{f'_{\lambda -1}}  & M_{\lambda } \ar[r] \ar[d]^{f'_{\lambda }} & M_{\lambda }/M_{\lambda -1} \ar[r] \ar[d]^{\bar{f'_{\lambda }}} & 0 \\
			0 \ar[r] & f(M_{\lambda -1})  \ar[r]  & f(M_{\lambda }) \ar[r]  & Im \ \bar{f'_{\lambda }}  \ar[r] & 0 
		}
	\end{equation}
	Since $\ker f'_{\lambda} = \ker f_{\lambda} \in \Filt{\X}$ and $ \ker \bar{f'_{\lambda }} = Im\ \alpha$ is $\X$-semisimple (by By Proposition  \ref{semisimple_modules} ), Proposition \ref{existence_X_socle_filtration} implies that $Im \ \bar{f'_{\lambda }}$ is also $\X$-semisimple. 
	\\%
	If $\lambda$ is a limit ordinal, then $f(M_{\lambda}) = \bigcup_{\alpha < \lambda} f(M_{\alpha})$. Hence 
	$Im \ f$ has an $Add(\X)$-filtration, proving that $\Filt{\X}$ is closed under extensions, kernels, and images. 
\end{proof}
The next proposition completes the proof that $\Filt{\X}$ is a wide subcategory of $\Mod{R}$.
\begin{proposition} \label{closed_under_cokernel}
		Let $\X$ be a semibrick, and let $M,N\in \Filt{\X}$ with $\X$-socle filtrations $(M_{\lambda}, \lambda \leq \mu)$ and $(N_{\omega}, \omega \leq \nu)$. Let $f\in \mod{R}(M,N)$ be a homomorphism.  
	\begin{itemize}
		\item[(a)] If two submodules of $N$ lie in $\Filt{\X}$, then both their sum and the intersection also belong to  $\Filt{\X}$.
		\item[(b)] The filtration $(Im\ f \cap N_{\lambda}, \lambda \leq \mu)$ is the $\X$-socle filtration of $Im\ f$. 
		\item[(c)] The cokernel of $f$ lies in $\Filt{\X}$.
	\end{itemize}
\end{proposition}
\begin{proof}
	(a) Let $U, V \subseteq N$ be submodules in $\Filt{\X}$. Consider the short exact sequence 
	$$
	\xymatrix{
		0 \ar[r] & U \cap V \ar[r] & U \oplus V \ar[r]^{\Sigma} & U+V\ar[r] & 0
	}
	$$
	where $\Sigma$ denotes the sum map. 
	Since $U \oplus V \in \Filt{\X}$ and $\Filt{\X}$ is closed under kernels and images (Proposition \ref{closed_under_kernel}), both $U+V$ and $U \cap V$ belong to $\Filt{\X}$.
	\\%
	(b) We show by transfinite induction on $\lambda$ that $(Im \ f \cap N_{\lambda}, \lambda \leq \mu)$ coincides with the $\X$-socle filtration $(V_{\lambda}, \lambda \leq \mu_0)$ of $Im \ f$. 
	Let $(V_{\lambda}, \lambda \leq \mu_0)$ be the $\X$-socle filtration of $f(M)$. For $\lambda = 1 $, we trivially have $soc_{\X} f(M) \subset f(M)\cap soc_{\X}\  N$. Conversely, $Im \ f \cap soc_{\X} N$ is in $\Filt{\X}$ and $\X$-semisimple by Proposition \ref{semisimple_modules}, hence equality holds. 
	\\%
	Assume the statement true for all ordinals less than $\lambda$. The inclusion $V_{\lambda} \subseteq f(M) \cap N_{\lambda}$  follows from Proposition \ref{strict-filtration}. For the converse, consider the commutative diagram
	\begin{equation} \label{higher-socles}
		\xymatrix{
			0 \ar[r] &  V_{\lambda -1} \ar[r] \ar[d]^{\varepsilon_{\lambda -1}} & f(M) \ar[r]^{\pi_{1}\hspace{0.3cm}} \ar[d]^{\varepsilon} & f(M) / V_{\lambda -1} \ar[r] \ar[d]^{\bar{\varepsilon}_{\lambda -1}} & 0 \\
			0 \ar[r] &  N_{\lambda -1} \ar[r] & N \ar[r]^{\pi_2\hspace{0.3cm}} & N/ N_{\lambda -1} \ar[r] & 0 
		}
	\end{equation}
		where $\varepsilon_{\lambda -1}$ denotes the inclusion. Since $V_{\lambda -1} = f(M) \cap N_{\lambda -1}$ by induction, $\bar{\varepsilon}_{\lambda -1}$ is injective. For any element $x\in Im \ f \cap N_{\lambda}$, its image $\pi_2(x)$ lies in $soc_{\X} (N/N_{\lambda -1})$; hence $\pi_1(x)$ lies in $soc_{\X} (N/N_{\lambda -1})$. This shows $Im \ f \cap N_{\lambda} \subseteq V_{\lambda}$. 
		\\%
		If $\lambda$ is a limit ordinal, then 
	$$f(M)\cap N_{\lambda}  =  f(M) \cap (\bigcup_{\alpha < \lambda} N_{\alpha}) = \bigcup_{\alpha < \lambda} f(M)  \cap N_{\alpha}= \bigcup_{\alpha < \lambda} V_{\alpha} = V_{\lambda} $$
	(c) Define a filtration $(W_\lambda , \lambda \leq \mu + \nu)$ of $N$ by 
	$$
	W_\lambda =
	\begin{cases}
		Im f \cap N_\lambda, & \text{if } \lambda \leq \mu,\\[4pt]
		Im f + N_\lambda, & \text{if } \mu < \lambda \leq \mu + \nu.
	\end{cases}
	$$
	We claim that $(W_\lambda , \lambda \leq \mu + \nu)$ is an $Add(\X)$-filtration, implying $\coker f \in \Filt{\X}$. For $\lambda \leq \mu$ this follows from (b). For $\mu <\lambda \leq \mu + \nu$, consider the diagram:
	$$ \xymatrix{
	0 \ar[r] & f(M) \cap N_{\lambda -1}  \ar[r] \ar[d] & f(M) \cap N_{\lambda} \ar[r] \ar[d] &  (f(M) \cap N_{\lambda}) /  (f(M) \cap N_{\lambda -1}) \ar[r] \ar[d] & 0 \\
	0 \ar[r]& f(M) \oplus N_{\lambda -1} \ar[r] \ar[d]^{\Sigma_{\lambda -1}} &  f(M) \oplus N_{\lambda } \ar[r] \ar[d]^{\Sigma_{\lambda }}  & (f(M) \oplus N_{\lambda }) /  (f(M) \oplus N_{\lambda -1}) \ar[r] \ar[d]^{\bar{\Sigma}_{\lambda -1}}  & 0 \\
	0\ar[r] & f(M)+ N_{\lambda -1} \ar[r] & f(M)+ N_{\lambda } \ar[r] &  (f(M)+ N_{\lambda }) /  (f(M)+ N_{\lambda -1}) \ar[r] & 0
	}$$
	By (b), $(f(M) \cap N_{\lambda}) /  (f(M) \cap N_{\lambda -1}) \in Add(\X)$. Since $ (f(M) \oplus N_{\lambda }) /  (f(M) \oplus N_{\lambda -1}) \cong N_{\lambda}/ N_{\lambda -1} \in Add(\X)$, Proposition \ref{semisimple} implies that  $(f(M)+ N_{\lambda } )/  (f(M)+ N_{\lambda -1})$ is in $Add(\X)$. Thus $(W_\lambda , \lambda \leq \mu + \nu)$ is an $Add(\X)$-filtration of $N$, proving $coker \ f \in \Filt{\X}$.
\end{proof}
The following Theorem summarises the Propositions \ref{closed_under_kernel} and \ref{closed_under_cokernel}:
\begin{theorem}
	For a semibrick $\X$, the category $\Filt{\X}$ is a wide subcategory of $\Mod{R}$. 
\end{theorem}
\section{$\X$-universal short exact sequences} \label{Universal_short_exact_sequences}
In this section we study the $X$-universal short exact sequence, which generalizes the construction introduced by Bongartz in \cite[Lemma 2.1]{Bongartz}  of his paper.
\begin{definition} \label{construction_universal_sequence}
	Let $\X$ be a semibrick and $Y\in \Mod{R}$. For each $X$ in the class $\X$, the extension group $\Ext^1(X,Y)$ carries a natural structure as a left $\End(X)$-module. Choose a basis of $\Ext^1(X,Y)$, say
	$$
	\xymatrix{
		(0 \ar[r] & Y \ar[r] & M_{X,i} \ar[r] & X \ar[r] & 0)_{i\in I_X}.
	}
	$$
	We then define the {\bf $\X$-universal short exact sequence starting from $Y$} as the pushout in the  diagram
	$$
	\xymatrix{
		0 \ar[r] & \bigoplus_{X\in \X} Y^{(I_X)} \ar[r] \ar[d]^{\sum} & \bigoplus_{X\in \X, i\in I_X} M_{X,i} \ar[r] \ar[d] & \bigoplus_{X\in \X} X^{(I_X)} \ar[r] \ar@{=}[d] & 0 \\
		0 \ar[r] & Y \ar[r]^{\alpha} & M \ar[r]^{\beta} & \bigoplus_{X\in \X} X^{(I_X)} \ar[r] & 0
	}
	$$
	An inclusion $\alpha: Y\hookrightarrow M$ is called an {\bf $\X$-universal inclusion} if the short exact sequence,
	$$ \xymatrix{
	0 \ar[r] &  Y \ar[r]^{\alpha} & M \ar[r]^{\pi \hspace{0.5cm}} & M/Im\ \alpha \ar[r] & 0
	}$$
	is an $\X$-universal short exact sequence starting from $Y$.
\end{definition}
\begin{lemma} \label{first_properties_X_universal_sequences}
	Let $\X$ be a semibrick and $Y\in \Mod{R}$. For a short exact sequence 
	\begin{equation}
			\label{short-exact-sequence}
		\xymatrix{
			0 \ar[r] & Y \ar[r]^{\alpha} & M  \ar[r]^{\beta\hspace{0.7cm}} & \bigoplus_{X\in \X} X^{(I_X)}  \ar[r] & 0
		}
	\end{equation}
	the following statements are equivalent:
	\begin{itemize}
		\item[(i)] The sequence is $\X$-universal.
		\item[(ii)] For every $X$ in $\X$, the connecting homomorphism $\delta : \Hom(X, \bigoplus_{X\in \X} X^{(I_X)} ) \ra \Ext^1(X, Y)$ is an isomorphism. 
	\end{itemize}
\end{lemma}
\begin{proof}
	$(i) \Rightarrow (ii)$: The surjectivity of $\delta$ follows directly from the construction. To show injectivity, let $f\in \Hom(X, \bigoplus_{X\in \X} X^{(I_X)} )$ be nonzero. Since $X$ is finitely generated, the image of $f$ is contained in a finite direct sum $\bigoplus_{i=1}^s X_i$. 
	For each $1\leq j \leq s$, let $\pi_j: \bigoplus_{i=1}^s X_i \ra X_j$ be the canonical projection and $\varepsilon_j:X_j\ra \bigoplus_{i=1}^s X_i$ the  inclusion. Then 
	$$ f= f\circ \sum_{i=1}^{s} \pi_i \circ \varepsilon_i = \sum_{i=1}^{s} (f\circ \pi_i) \circ \varepsilon_i = \sum_{i=1}^{s} \lambda_i \cdot \varepsilon_i$$
	with $\lambda_i=f\circ \pi_i \in \End(X)$. Applying $\delta$ yields $\delta (f) = \sum_{i=1}^{s} \lambda_i \cdot \delta(\varepsilon_i)$. The elements $(\delta(\varepsilon_i))_{1\leq i \leq s}$ form a linearly independent family in $\Ext^1(X,Y)$,  so $\delta(f)=0$ implies $f=0$.
	\\%
	$(ii) \Rightarrow (i)$:  For each $X$ in $\X$ and each $i\in I_X$, let $\varepsilon_{X,i}: X \ra X^{(I_X)}$ denote the inclusion and $E_{X,i} = \delta (\varepsilon_{X,i})$. Since $(\varepsilon_{X,i})_{i\in I_X}$ is a basis of $\Hom(X, X^{(I_X)} )$  as an $\End(X)$-module, and $\delta$ is an isomorphism, the elements $(E_{X,i})_{i\in I_X}$  form an $\End(X)$-basis of $\Ext^1(X,Y)$. Using this basis in the construction of  the $\X$-universal short exact sequence yields exactly the given sequence, completing the proof. 
\end{proof}
Happel and Unger observed in \cite[§1]{happel1989almost} that the universal sequence associated with $\Ext^1(X,Y)$ is uniquely determined up to isomorphism whenever $X$ is a brick. The next result shows, that this uniqueness also holds for the generalized $\X$-universal short exact sequence. 
\begin{proposition} \label{unique_bongartz}
	Let $\X$ be a semibrick and an $R$-module. Then all $\X$-universal short exact sequences starting from $Y$ are isomorphic. 	
\end{proposition}
\begin{proof}
	Let $E$ and $E'$ be two $\X$-universal short exact sequences:
	$$
	\xymatrix@C=1pc @R=.5pc{
		E: \hspace{0.3cm}  0\ar[r] &  Y  \ar[r]^{\alpha}  & Z\ar[r]^{\beta\hspace{0.8cm}} & \bigoplus_{X\in \X} X^{(I_X)}  \ar[r] & 0\\
		E': \hspace{0.3cm}  0\ar[r] &  Y  \ar[r]^{\alpha'}  & Z' \ar[r]^{\beta'\hspace{0.8cm}} & \bigoplus_{X\in \X} X^{(I_X)}  \ar[r] & 0.
	}
	$$
	We have to show that there exist isomorphisms $f:Z' \ra Z$ and $g\in \End(\bigoplus_{X\in \X} X^{(I_X)})$ making the diagram commute. 
	Let 
	$$ (E_{X,i}: \hspace{0.3 cm} 0 \lra Y  \stackrel{\alpha_{X,i}}{\lra} Z_{X,i} \stackrel{\beta_{X,i}}{\lra} X \lra 0)_{i \in I_X, X\in \X}$$
	and 
	$$ (E_{X,i}': \hspace{0.3 cm} 0 \lra Y  \stackrel{\alpha_{X,i}'}{\lra} Z_{X,i} '\stackrel{\beta_{X,i}'}{\lra} X \lra 0)_{i \in I_X, X\in \X}$$
	be the $\End(X)$-bases of $\Ext^1(X, Y)$ used in the constructions of $E$ and $E'$. 
	\\%
	For $X\in \X$, apply $\Hom(-,X)$ to the sequence $E$. Since the connecting homomorphism $\delta$ is surjective, for every $i\in I_X$ there exists a homomorphism $g_{X,i}\in \Hom(X,\bigoplus_{X\in \X} X^{(I_X)} )$ satisfying $\delta (g_{X,i}) = E_{X,i}'$.
	
	Let $g_X\in \End(X^{(I_X)})$ be the induced endomorphism such that $\varepsilon_i \circ g_X= g_{X,i}$ 
	for all $i\in I_X$, where $\varepsilon_i : X \ra X^{(I_X)}$ are the canonical inclusions. 
	\\%
	We first claim that $g_X$ is surjective. Each $\varepsilon_i$ can be expressed as a finite linear combination of the homomorphisms $g_{X,j}$. Since $(E_{X,i}')_{i\in I_X}$ is an $\End(X)$-basis of $\Ext^1(X,Y)$, every $E_{X,i}$ can be written as $ \sum_{j=1}^{r} \lambda_j \cdot E_{X,j}'$ for some  $\lambda_j\in \End(X)$. Then 
	$$\delta( \varepsilon_i)= E_{X,i} = \sum_{j=1}^{r} \lambda_j \cdot E_{X,j}' =  \sum_{j=1}^{r} \lambda_j \cdot \delta (g_{X,j})=  \delta (\sum_{j=1}^{r} \lambda_j  \cdot g_{X,j}).$$
	As $\delta$ is an isomorphism (by  Lemma \ref{first_properties_X_universal_sequences}(ii)), we obtain $\varepsilon_i= \sum_{j=1}^{r} \lambda_j  g_{X,j}$, proving that $g_X$ is surjective. 
	\\%
	Now we show that $g_X$ is injective. Take a nonzero map $f: X \ra X^{(I_X)}$ with  $Im\ f \subset soc_{\X}\ (ker (g_X))$. Since the image of $f$ lies in a finite direct sum, say $Im\ f \subset \bigoplus_{j=1}^r X$ we may write $ f =  \sum_{j=1}^r \lambda_j \cdot  \varepsilon_j$ with $\lambda_j\in \End(X)$ not all zero. Then  
	$$ 0 = f\circ g_X = (\sum_{j=1}^r \lambda_j \cdot  \varepsilon_j) \circ g_X = \sum_{j=1}^r \lambda_j \cdot ( \varepsilon_j \circ g_X ) = \sum_{j=1}^r \lambda_j \cdot  g_{X,i},$$
	contradicting the linear independence of the homomorphisms $(g_{X,i})_{i\in I_X}$. 
	Hence $g_X$ is an isomorphism. Consequently, the direct sum $g= \bigoplus_{X\in \X} g_X $ is also an isomorphism, and the pull-back of $E$ along $g$ yields $E'$. 
\end{proof}
\begin{definition} \label{Y_X_infty}
	Let $\X$ be a semibrick and $Y\in \Mod{R}$. Define a direct system $(Y_{\X}(i), \alpha_{i, i+1} \ |\  i\in \mathbb{N})$ recursively by setting $Y_{\X}(1)=Y$, and for each $i$, let $\alpha_{i,i+1} : Y_{\X}(i)\ra Y_{\X}(i+1) $ be the monomorphism of an $\X$-universal sequence starting from $Y_{\X}(i)$. The direct limit of this system is denoted by $Y_{\X}(\infty) $. 
\end{definition}
Given two such direct systems $(Y_{\X}(i), \alpha_{i, i+1} \ |\  i\in \mathbb{N})$ and $(Y'_{\X}(i), \alpha'_{i, i+1} \ |\  i\in \mathbb{N})$ with the same starting module $Y_{\X}(1)= Y_{\X}'(1)$,  Proposition \ref{unique_bongartz} ensures the existence of isomorphisms  $f_i: Y_{\X}(i) \ra Y_{\X}(i)$ for $i\geq 2$ that make the diagram 
$$\xymatrix{
Y \ar[r]^{\alpha_{1,2} \hspace{0.4cm}} \ar@{=}[d] & Y_{\X}(2) \ar[r]^{\alpha_{2,3} } \ar[d]^{f_2} & Y_{\X}(3) \ar[r] \ar[d]^{f_3} & \hdots \\
Y \ar[r]^{\alpha'_{1,2} \hspace{0.4cm}} & Y_{\X}(2) \ar[r]^{\alpha'_{2,3}} & Y_{\X}(3) \ar[r]& \hdots
}$$
commute. Thus the direct limit $\alpha: Y \ra Y_{\X}(\infty)$ is well-defined up to isomorphism, and we may view $Y$ as a submodule of $Y_{\X}(\infty)$.
\begin{proposition} \label{Hom-property}
	Let $\X$ be a semibrick and $Y\in \Mod{R}$. For any $f\in \Hom(X, Y_{\X}(\infty))$ with $X\in \X$, the image of $f$ is contained in $Y$.
\end{proposition}
\begin{proof}
	For $r\in \mathbb{N}$, consider the  $\X$-universal short exact sequence
	$$ \xymatrix{
		0 \ar[r] & Y_{\X}(r) \ar[r]^{\alpha\hspace{0.5cm}} & Y_{\X} (r+1) \ar[r]^{\beta\hspace{0.7cm}} & Y_{\X} (r+1) / Y_{\X}(r) \ar[r] & 0.
	} $$
	If there existed a homomorphism $f: X\ra Y_{\X} (r+1)$ with $f\circ \beta \ne 0$, then the connecting homomorphism $\delta: \Hom(X, Y_{\X} (r+1) / Y_{\X}(r)) \ra \Ext^1(X, Y_{\X}(r) )$ would satisfy $\delta (f \circ \beta)=0$, contradicting Lemma  \ref{first_properties_X_universal_sequences}. Hence $Im \ f \subset Y$.
\end{proof}
From now on we need to assume that each brick in the semibrick $\X$ has the projective dimension at most one. 
\begin{lemma}  \label{push-out_along_universal_inclusion}
Let $\X$ be a semibrick with $\pd{\X} \leq 1$. Consider the commutative diagram with exact rows and columns 
\begin{equation} \label{cross-universal-sequences}
	\xymatrix{
		&& (i) & (iii) \\
		& & 0 \ar[d] & 0 \ar[d]\\
		&	& U \ar[r]^{id} \ar[d]^{f}& U  \ar[d]^{f'} \\
		(ii) &	0 \ar[r]& Y \ar[r]^{\alpha} \ar[d]^{g} & Y_{\X}(2) \ar[r]^{\beta} \ar[d]^{g'} & Y_{\X}(2) / Y \ar[r] \ar@{=}[d] & 0 \\
	(iv)	&	0 \ar[r] & Y/ U \ar[r]^{\alpha'} \ar[d] & Y_{\X}(2) / U \ar[r]^{\beta'}  \ar[d] & Y_{\X}(2) / Y \ar[r]  & 0 \\
	&& 0 & 0
	}
\end{equation} 
Assume that the short exact sequence $(ii)$ is $\X$-universal, and that in sequence $(i)$ the connecting homomorphism $\delta_1 : Hom(X,Y/ U ) \ra \Ext^1(X, U)$ is surjective for every $X\in \X$.
\\%
Then the connecting homomorphism $\Hom(X, Y_{\X}(2) / U ) \ra \Ext^1(X, U)$ from $(iii)$ is also surjective for all $X\in \X$, and the short exact sequence $(iv)$ is $\X$-universal.
\end{lemma}
\begin{proof}
For $X\in \X$ we apply the functor $\Hom(X,-)$ to the diagram \eqref{cross-universal-sequences}. This gives a diagram of long exact sequences
$$ \xymatrix{
	& (X , U) \ar[r]^{id} \ar[d]& (X , U)  \ar[d]  & & {}^1 (X,U) \ar[d]^{ {}^1(X,f) }\\
	0 \ar[r]& (X , Y) \ar[r] \ar[d] & (X , Y_{\X}(2)) \ar[r]^{(X,\beta)} \ar[d]^{(X,g')} & (X , Y_{\X}(2) / Y) \ar[r]^{\delta_3} \ar@{=}[d] & {}^1(X , Y) \ar[d]^{{}^1(X, g)}  \\
	0 \ar[r] & (X , Y/ U) \ar[r] \ar[d]^{\delta_1} & (X , Y_{\X}(2) / U) \ar[r]^{(X, \beta')}  \ar[d]^{\delta_2} & (X , Y_{\X}(2) / Y) \ar[r]^{\delta_4} &  {}^1(X , Y / U) \ar[d]  \\
	0 \ar[r]& {}^1(X, U) \ar[r]^{\cong} & {}^1(X, U) &   & {}^2(X, U) 
 } $$
 It is assumed that the connecting homomorphism  $\delta_1: \Hom(X,Y/U) \ra \Ext^1(X,U)$ is surjective. This implies that $\delta_2$ is also surjective; in particular, the connecting homomorphism $\Hom(X, Y_{\X}(2) / U ) \ra \Ext^1(X, U)$ in (iii) is surjective for all $X$ in $\X$.
 \\%
 The connecting homomorphism $\delta_3$ is an isomorphism because the short exact sequence (ii) is $\X$-universal. Since $\pd{\X}\leq 1$, we have $\Ext^2(X, U)=0$. As $\delta_1$  is surjective, $\Ext^1(X,f)$ must be the zero map, which implies that $\Ext^1(X,g)$ is an isomorphism. Consequently, $\delta_4$ is an isomorphism, and by Lemma \ref{first_properties_X_universal_sequences} the short exact sequence $(iv)$ is $\X$-universal. 
 \end{proof}
\begin{proposition} \label{quotient_universal_sequences}
Let $\X$ be a semibrick with $\pd{\X} \leq 1$.
\begin{itemize}
	\item[(a)] Let $Y$ be an $R$-module and $Y\hookrightarrow Y_{\X}(2) \hookrightarrow Y_{\X}(3) \hookrightarrow \hdots \ $ be a chain of $\X$-universal inclusions. Then, for each $j\geq 1$, the induced chain $Y_{\X}(j+1) / Y_{\X}(j) \hookrightarrow Y_{\X}(j+2) / Y_{\X}(j)  \hookrightarrow \hdots \ $ is also a chain of $\X$-universal inclusions.
	 \item[(b)] Let $Y\in Add(\X)$ and $r\in \mathbb{N}$. Then the modules $(Y_{\X}(i)$ for $i \leq r)$ form the $\X$-socle filtration of $Y_{\X}(r)$. 
\end{itemize}
\end{proposition}
\begin{proof}
	(a) For $j=1$ we apply  Lemma \ref{push-out_along_universal_inclusion} to  the two short exact sequences 
$$ \xymatrix@C=1pc @R=.5pc{
(i) & 0 \ar[r] & Y \ar[r] & Y_{\X}(2) \ar[r] & Y_{\X}(2) /Y \ar[r] & 0 \\
(ii) & 0 \ar[r] & Y_{\X}(2) \ar[r] & Y_{\X}(3) \ar[r] & Y_{\X}(3) /Y_{\X}(2) \ar[r] & 0
}$$
This yields two further sequences:
$$ \xymatrix@C=1pc @R=.5pc{
	(iii) & 0 \ar[r] & Y \ar[r] & Y_{\X}(3) \ar[r] & Y_{\X}(3) /Y \ar[r] & 0 \\
	(iv) & 0 \ar[r] & Y_{\X}(2)/Y \ar[r] & Y_{\X}(3)/Y \ar[r] & Y_{\X}(3) /Y_{\X}(2) \ar[r] & 0
}$$
The connecting homomorphism of sequence $(iii)$ is surjective for all $X\in \X$, and $(iv)$ is $\X$-universal. We can use Lemma  \ref{push-out_along_universal_inclusion} for the two short exact sequences $(i), (ii)$:
$$ \xymatrix@C=1pc @R=.5pc{
	(i')  & 0 \ar[r] & Y \ar[r] & Y_{\X}(3) \ar[r] & Y_{\X}(3) /Y \ar[r] & 0  \\
	(ii') & 0 \ar[r] & Y_{\X}(3) \ar[r] & Y_{\X}(4) \ar[r] & Y_{\X}(4) /Y_{\X}(3) \ar[r] & 0
}$$
and get the next step. Iterating the construction proves the claim for $j=1$. For $j=2$, apply the same reasoning to the factor module $Y_{\X}(3)/ Y_{\X}(2) \cong (Y_{\X}(3)/Y) / (Y_{\X}(2) / Y)$ and similarly for higher $j$. 
\\%
(b) We proceed by induction on $i$. For $i=1$, Proposition  \ref{Hom-property} implies $soc_{\X} Y_{\X}(r) = Y$. Assume the statement holds up to $i$. Consider the short exact sequence 
$$ \xymatrix{
	0 \ar[r] & Y_{\X}(i) \ar[r] & Y_{\X}(r) \ar[r] & Y_{\X}(r) / Y_{\X}(i) \ar[r] & 0
}$$ 
By part (a) the chain 
$$ Y_{\X}(i+1) / Y_{\X}(i) \hookrightarrow Y_{\X}(i+2) / Y_{\X}(i) \hookrightarrow \hdots Y_{\X}(r) / Y_{\X}(i) $$ 
consists of $\X$-universal inclusions. Hence, for any $X$ in $\X$ and any homomorphism $f:X\ra Y_{\X}(r) / Y_{\X}(i)$, the image of $f$ lies in $Y_{\X}(i+1) / Y_{\X}(i)$. Therefore $soc_{\X} \ (Y_{\X}(r) / Y_{\X}(i)) = Y_{\X}(i+1) / Y_{\X}(i)$ showing that $Y_{\X}(i+1)$ is the next submodule in the $\X$-socle filtration. 
\end{proof}
\begin{corollary} \label{homomorphism_of_filtered_objects}
Let $\X$ be a semibrick with $\pd{\X} \leq 1$, and let $U, V$ be $R$-modules. Let $f: U_{\X}(\infty) \ra V_{\X}(\infty)$ be a homomorphism such that $f(U) \subset V_{\X}(r)$ for some $r\in \mathbb{N}$. Then $f(U_{\X}(i)) \subset V_{\X}(i+r-1)$ for all $i\in \mathbb{N}$.  
\\%
In particular, if  $U$ and  $V$ lie in $\Filt{\X}$ and  $U$ is $\X$-semisimple, then $f(U_{\X}(i)) \subset V_{\X}(i)$ for all $i\in \mathbb{N}$.
\end{corollary}
\begin{proof}
	Let $\pi_1: U_{\X}(\infty) \ra U_{\X}(\infty)/U$ and $\pi_2: V_{\X}(\infty) \ra V_{\X}(\infty)  / V_{\X}(r)$ be the canonical projections. Because  $f(U) \subset V_{\X}(r)$, the map $f$ induces a homomorphism $f\circ \pi_2$ over $U_{\X}(\infty)/U$ making the corresponding diagram commute.
$$ \xymatrix{
U_{\X}(\infty) \ar[r]^f \ar[d]^{\pi_1} & V_{\X}(\infty) \ar[d]^{\pi_2} \\
U_{\X}(\infty)/U \ar[r]^{\bar{f}} & V_{\X}(\infty) / V_{\X}(r)
}$$
By Proposition  \ref{quotient_universal_sequences}(b), the modules on both sides carry $\X$-socle filtrations compatible with $\bar{f}$, meaning $\bar{f}$ is a homomorphism of $\X$-filtered objects. Thus 
$f(U_{\X}(i)) \subset V_{\X}(i+r-1)$ for all $i\in \mathbb{N}$. 
\\%
If, in addition, $U$ and $V$ lie in $\Filt{\X}$ and $U$ is $\X$-semisimple, then $f(U_{\X}(1)) \subset V_{\X}(1)$, which gives the stated special case. 
\end{proof}
\section{$\X^{\perp}$-envelopes of modules for a semibrick $\X$}  \label{Envelopes_of_modules_for_a_semibrick}
For a class $\mathcal{C}$ of modules, the following notions of envelopes generalize the classical concept of injective modules. We recall the definitions of a $\mathcal{C}$-envelope from \cite[Definition 5.1]{gobel2012approximations}. 
\begin{definition}
	Let $\mathcal{C}$ be a class of modules. A homomorphism $f\in \Hom(Y,C)$ with $C\in \mathcal{C}$ is called a {\bf $\mathcal{C}$-preenvelope of $Y$}, if for every homomorphism $f': Y \ra C'$ for a $C'\in \mathcal{C}$, there exists a homomorphism $g:C\ra C'$ such that $f'=f\circ g$:
	$$
	\xymatrix{
		Y \ar[r]^{f}  \ar[dr]_{f'}  & C \ar@{-->}[d]^{g} \\
		& C'	
	}
	$$
	A $\mathcal{C}$-preenvelope $f$ is called a {\bf special $\mathcal{C}$-preenvelope}, if $f$ is injective and the cokernel of $f$ lies in the class ${}^{\perp}{\mathcal{C}}$, that is, $\Ext^1(\coker f, C)=0$ for all $C\in \mathcal{C}$. 
	\\%
	A $\mathcal{C}$-preenvelope $f$ is a {\bf $\mathcal{C}$-envelope of $Y$} if it is left minimal; that is, for every endomorphism $g\in \End(C)$, the equality $f = f\circ g$ implies that $g$ is an automorphism. 
\end{definition}
It is well known \cite[Theorem 6.11]{gobel2012approximations} that for every module $M\in \Mod{R}$ and any set of modules $\mathcal{S}$ there exists a short exact sequence 
$$ \xymatrix{
0 \ar[r] & M \ar[r] & P \ar[r] & N \ar[r] & 0
} $$
with $P\in \mathcal{S}^{\perp}$ and $N\in \Filt{\mathcal{S}}$. In particular, the embedding  $M \ra P$ is a special $\mathcal{S}^{\perp}$-preenvelope $M$.  
\begin{theorem} \label{envelope}
	Let $\X$ be a semibrick with $\pd{\X} \leq 1$. Then for $Y\in \Mod{R}$,  the embedding $\alpha:Y \ra Y_{\X}(\infty)$ is an $\X^{\perp}$-envelope of $Y$. 
\end{theorem}
\begin{proof}
First we show that $\alpha: Y\hookrightarrow Y_{\X}(\infty)$ is an $\X^{\perp}$-preenvelope of $Y$. Let $f:Y\ra C$ be a homomorphism with  $C\in \X^{\perp}$. Applying the functor $\Hom(-, C)$ to the  short exact sequence
$$
\xymatrix{
	0 \ar[r] & Y \ar[r]^{\alpha \hspace{0.5cm}}  & Y_{\X}(\infty) \ar[r]^{\beta \hspace{0.4cm}} & Y_{\X}(\infty)/Y \ar[r] & 0 
}
$$
gives part of a long exact sequence 
$$ \hdots  \lra (Y_{\X}(\infty),C) \lra (Y,C) \lra {}^1(Y_{\X}(\infty)/Y , C) \lra \hdots $$
Since $Y_{\X}(\infty)/Y$ lies in $\Filt{\X}$ and $C$ lies in $\X^{\perp}$, the Eklof Lemma \ref{eklof-lemma} implies $\Ext^1(Y_{\X}(\infty)/Y , C)=0$. Hence there exists $g\in \Hom(Y_{\X}(\infty), C)$ with $f = \alpha \circ g$.
\\%
Next, we show that $Y_{\X}(\infty)$ itself lies in $X^{\perp}$. For each $i\in \mathbb{N}$, consider the $\X$-universal short exact sequence:
$$\xymatrix{
0 \ar[r] & Y_{\X}(i-1) \ar[r]^{\alpha_i} & Y_{\X}(i) \ar[r]^{\beta_i \hspace{1cm} } & Y_{\X}(i) / Y_{\X}(i-1) \ar[r] & 0
}$$
By Lemma \ref{first_properties_X_universal_sequences}(ii) the connecting homomorphism $\delta_i : \Hom(X, Y_{\X}(i)/Y_{\X}(i-1)) \ra \Ext^1(X, Y_{\X}(i-1))$ is an isomorphism for all $X\in \X$. This implies that $\Ext^1(X, \alpha_i)$ maps trivially into $\Ext^1(X, Y_{\X}(i))$. By Lemma \ref{direct_limes_of_Ext}, we obtain 
$$ \Ext^1(X, Y_{\X}(\infty)) = \Ext^1(X, \lim_{i\ra \infty} Y_{\X}(i)) = \lim_{i\ra \infty} \Ext^1(X, Y_{\X}(i))=0.$$
Thus $Y_{\X}(\infty)$ lies in $X^{\perp}$.
\\%
Now we prove that $\alpha$ is an $\X^{\perp}$-envelope. Assume $\alpha = \alpha \circ h$ for some endomorphism $h\in \End(Y_{\X}(\infty))$. Applying the Snake Lemma to the diagram
$$
\xymatrix{
	&&\ker h \ar[d] \ar[r]^{\cong} & \ker  \bar{h} \ar[d] \\
	0 \ar[r] & Y \ar[r]^{\alpha \hspace{0.5cm} }  \ar@{=}[d] & Y_{\X}(\infty) \ar[r]^{\beta} \ar[d]^{h}& Y_{\X}(\infty)/Y \ar[r] \ar[d]^{\bar{h}}& 0 \\
	0 \ar[r] & Y \ar[r]^{\alpha \hspace{0.5cm}} & Y_{\X}(\infty) \ar[r]^{\beta} \ar[d] & Y_{\X}(\infty)/Y \ar[r] \ar[d] & 0 \\
	&& \coker h  \ar[r]^{\cong} & \coker  \bar{h} 
}
$$
shows that $\ker  h \cong \ker  \bar{h}$ and $\coker h \cong \coker  \bar{h}$. 
By Propositions \ref{closed_under_kernel},  $\ker \bar{h}$ lies in $\Filt{\X}$, and hence so does $\ker h$. From Proposition \ref{Hom-property}, we know that $soc_{\X} (\ker  h)$ is contained in $\alpha (Y)$; this forces $\beta (soc_{\X} (\ker h)) =0$, which is possible only if $\ker  \bar{h} =0$. Thus $h$ is injective. 
\\%
Similarly, Proposition \ref{closed_under_cokernel} implies that $\coker h$ lies in $\Filt{\X}$. Therefore $\Ext^1( \coker h , Y_{\X}(\infty))=0$, and there exists a submodule $U \subset Y_{\X}(\infty)$ such that $Y_{\X}(\infty)= U \oplus h(Y_{\X}(\infty))$. Since $U\cong \coker h \in \Filt{\X}$,  Proposition \ref{Hom-property} implies that for all $X\in \X$ and $f \in \Hom(X, Y_{\X}(\infty))$ we have $Im \ f \subset \alpha (Y) \subset h(Y_{\X}(\infty))$. This can happen only if $U=0$, so $h$ is surjective. Hence $h$ is an automorphism, proving that $\alpha$ is left minimal. Therefore $\alpha$ is an $\X^{\perp}$-envelope of $Y$. 
\end{proof}
\begin{remark}
	Following an observation by Lidia Angeleri H\"ugel: Let $\X$ be a semibrick with $\pd{\X} \leq 1$ and $\Hom(\X , R) =0$. Let $\X^{\perp_{0,1}}$ denote the class of modules $M\in \Mod{R}$ satisfying  $Hom(\X,M)=\Ext^1(\X,M)=0$. 
	\\%
	For any module $U\in \X^{\perp_{0,1}}$, applying the functor $\Hom(-,U)$ to the short exact sequence 
		$$
	\xymatrix{
		0 \ar[r] & R \ar[r]^{\alpha \hspace{0.5cm}} & R_{\X}(\infty) \ar[r]^{\beta \hspace{0.3cm}} & R_{\X}(\infty)/R \ar[r] & 0
	}
	$$
	yields an isomorphism $\Hom (R, U ) \cong \Hom(R_{\X}(\infty), U)$. Thus the embedding 
	$R\ra R_{\X}(\infty)$ is a $\X^{\perp_{0,1}}$-reflection, in the sense of  \cite[page 712]{hugel2019characterisation}. By \cite[Proposition 2.7(1)]{hugel2019characterisation} in that paper, the category $\X^{\perp_{0,1}}$ is bireflective. Hence, by \cite[Theorem 2.4]{hugel2019characterisation}, the image of the restriction functor associated with the ring epimorphism $R\ra R_{\X}(\infty)$ coincides with the reflective subcategory $\X^{\perp_{0,1}}$.
	Since the reflection is unique, $R\ra R_{\X}(\infty)$ must be a ring epimorphism realizing this reflection. Finally, \cite[Proposition 2.7(2)]{hugel2019characterisation} states that this map is the universal localization of $R$ at $\X$. 
	\\%
	The explicit construction of $R_{\X}(\infty)$ therefore provides a concrete method to analyze universal localizations of rings at a given semibrick $\X$. 
\end{remark}
\section{Pr\"ufer modules} \label{Pruefer_modules}
The notion of a Pr\"ufer module was first introduced by Claus Michael Ringel in 1979 in his paper on infinite-dimensional modules over hereditary algebras (\cite[p. 326]{ringel1979infinite}). He investigated direct limits of irreducible monomorphisms inside tubes, showed that they are indecomposable and $\Sigma$-pure-injective, and called them Prüfer modules, in analogy with Pr\"ufer groups.
Later, inspired by work of Henning Krause, Ringel generalized the concept: a nonzero module $M$ over an arbitrary ring is called a Pr\"ufer module if it possesses a locally nilpotent, surjective endomorphism whose kernel has finite length. Such modules need not be indecomposable, but they remain $\Sigma$-pure-injective. Moreover, for a finite-dimensional algebra, the product closure of a Pr\"ufer module of infinite type always contains a generic module (\cite{ringel2019construction}).
\\%
We now propose a homological generalization of the notion of Pr\"ufer modules.
To avoid calling finitely generated modules Pr\"ufer modules, we assume that the semibrick $\X$  contains no $\X$-divisible objects. Here, a module $M$ is called $\X$-divisible if $\Ext^1(\X,M)=0$.
Under this assumption, for each $Y\in \X$ and every $r\in \mathbb{N}$ there exists some $X\in \X$ such that $\Ext^1(X, Y_{\X}(r)) \ne 0$. Hence, $Y_{\X}(\infty)$ is  not finitely generated.
\begin{definition}
		Let $\X$ be a semibrick with projective dimension $\pd{\X} \leq 1$ that contains no $\X$-divisible objects. We call the $\X^{\perp}$-envelope $Y_{\X}(\infty)$ of a brick $Y\in \X$ a {\bf Pr\"ufer module}. If there is no risk of confusion about the underlying semibrick $\X$, we write $\bar{Y}$ instead of $Y_{\X}(\infty)$.
\end{definition}
\begin{lemma} \label{Pruefer-modules_have_simple_socle}
	Let $\X$ be a semibrick with $\pd{\X} \leq 1$ and without $\X$-divisible objects. Then for every $Y\in \X$ we have $soc_{\X} \ \bar{Y} =Y$. In particular, every Pr\"ufer module is indecomposable. 
\end{lemma}
\begin{proof}
	The simplicity of the $\X$-socle of $Y_{\X} (\infty)$ follows directly from Proposition \ref{Hom-property}, which states that the image of any homomorphism from an object in $\X$ to $Y_{\X} (\infty)$ is contained in $Y$. 
\end{proof}
\begin{proposition}
	Let $\X$ be a semibrick with $\pd{\X} \leq 1$ and without $\X$-divisible objects, and let $Y\in \X$. Then the endomorphism ring $E=\End(\bar{Y})$ of the Pr\"ufer module $\bar{Y}$ is local, and its  Jacobson  radical is $J =\{ f\in E\  | \ f(Y)=0 \}$. Moreover, every nonzero endomorphism of $\bar{Y}$ is surjective. 
\end{proposition}
\begin{proof}
	We show that $J$ contains all non-invertible elements of $E$. Since the sum of two such endomorphisms again sends $Y$ to zero, $J$ is an ideal. Thus $E$ is a local ring with radical $J$.
	\\%
	Let $f\in E$ be an endomorphism such that $f(Y)\ne 0$. Because $Y= soc_{\X} \bar{Y}$ and $f(soc_{\X} \bar{Y}) \subset soc_{\X} \bar{Y}$, we have  $f(Y) \subset Y$, and since $Y$ is a brick,  $f(Y)=Y$. Consider the following exact diagram of $f$ and its restriction $f'$ to $Y$:
	$$ \xymatrix{
		&  & \ker f \ar[r]^{\cong} \ar[d] & \ker \bar{f}  \ar[d] &  \\ 
		0 \ar[r] & Y \ar[r]  \ar[d]^{f'} & \bar{Y} \ar[r]^{\beta} \ar[d]^f &  \bar{Y}/Y \ar[r] \ar[d]^{\bar{f}} & 0 \\
		0 \ar[r] & Y \ar[r] & \bar{Y} \ar[r] \ar[d] &  \bar{Y}/Y \ar[r] \ar[d] & 0\\
		&  & \coker f \ar[r]^{\cong} & \coker \bar{f}  & \\ 
	}$$
	Both $\ker \bar{f}$ and $\coker \bar{f}$ lie in $\Filt{\X}$; hence $\ker f$ and $\coker f$ also lie in $\Filt{\X}$. The $\X$-socle of $\ker f$ maps into $Y$, and the map $\beta$ restricted to $soc_{\X} \ (\ker f)$ must be injective. This is only possible if $\ker f =0$, so $f$ is injective. The cokernel of $f$ lies in $\Filt{\X}$, implying that $f$ splits. Since $ \bar{Y} $ is indecomposable, this can happen only if $\coker f=0$, so $f$ is surjective. Therefore $J$ consists precisely of the noninvertible endomorphisms. The last argument also showsthat every nonzero endomorphism of $\bar{Y}$ is surjective. 
\end{proof}
\begin{remark}
Let $\X$ be a semibrick that contains no $Y$ with $\Ext^1(Y,\X)=0$. Then $\Filt{\X}$ has no nonzero   $\Ext$-projective objects. Indeed, suppose there exists a nonzero module $M\in \Filt{\X}$ that is $\Ext$-projective in $\Filt{\X}$. For $X\in \X$ consider the short exact sequence 
$$ \xymatrix{
	0 \ar[r] & soc_{\X} M \ar[r] & M \ar[r] & M/soc_{\X} M \ar[r] & 0.
}$$
Applying the functor $\Hom(-,X)$ yields a segment of the long exact sequence
$$ \xymatrix{ \hdots \ar[r] & \Ext^1(M,X) \ar[r] & \Ext^1(soc_{\X} M ,X) \ar[r] \ar[r]& 0 },$$
since $\Ext^2(M/soc_{\X}M) , X) =0$ (because  $\pd{(M/soc_{\X}M)} \leq 1$ ). By assumption, there exists some $X\in \X$ with $\Ext^1(soc_{\X} M, X)\ne 0$,  which contradicts the $\Ext$-projectivity of $M$. Therefore, $\Filt{\X}$ has no nonzero $\Ext$-projective modules. 
\end{remark}
The next theorem provides a classification of the $\Ext$-injective objects in $\Filt{\X}$.
\begin{theorem} \label{classification_Pruefer_modules}
Let $\X$ be a semibrick with $\pd{\X} \leq 1$ and without $\X$-divisible objects. Let $M$ be a module in $\Filt{\X}$ whose $\X$-socle decomposes as $soc_{\X}M=\oplus_{i\in I} M_i$, where each $M_i$ is isomorphic to a brick in $\X$. 
\begin{itemize}
	\item[(a)]  For each $i\in I$, let $\alpha_i:M_i\ra \bar{M_i}$ denote the $X^{\perp}$-envelope of $M_i$. Then there exists a monomorphism $f$ making the following diagram commute:
	$$
	\xymatrix{
		\oplus_{i\in I} M_i \ar@{^(->}[r]^{ \hspace{0.4cm}\varepsilon}     \ar[d]_{\oplus_{i\in I} \alpha_i} & M \ar@{-->}[ld]^{f \vspace{0.7cm}} \\
		\oplus_{i\in I} \bar{M_i}&
	}
	$$
The homomorphism $f: M \ra \oplus_{i\in I} \bar{M_i}  $ is an $\X^{\perp}$-envelope of $M$.
	\item[(b)] If $M$ is $\Ext$-injective in $\Filt{\X}$, then $M$ is isomorphic to a direct sum of Pr\"ufer modules.  
\end{itemize}
\end{theorem}
\begin{proof}
	(a) Since $\oplus_{i\in I} \bar{M_i}$ is $\Ext$-injective in $\Filt{\X}$, the homomorphism $\oplus_{i\in I} \alpha_i$ can be lifted to a homomorphism $f$. Because $\coker f$ lies in $\Filt{\X} \subset {}^{\perp}(\X^{\perp })$, it remains to show that $f$ is injective and left minimal.
	By construction, the restriction of $soc_{\X} M$ is injective. To show left minimality of $f$, suppose there  exists a $g\in \End(\oplus_{i\in I} \bar{M_i})$ such that $f= f\circ g$.  
	\\%
	Since $f$ is injective, $g$ is injective on $f(M)$, in particular on the submodule $\oplus_{i\in I} M_i= soc_{\X} (\oplus_{i\in I} \bar{M_i}) $. As $soc_{\X} (\ker g) \subset soc_{\X} (\oplus_{i\in I} \bar{M_i}) $, we have $soc_{\X } (\ker g)=0 $, so $g$ is injective. 
	\\%
	Since $\coker g$ lies in $\Filt{\X}$, the monomorphism $g$ splits. Write $\oplus_{i\in I} \bar{M_i} = Im\ g \oplus C$, where $C$ is a submodule in $\Filt{\X}$. By construction, $soc_{\X} (\oplus_{i\in I} \bar{M_i} ) = f(\oplus_{i\in I} M_i ) \subset Im \ g$. Hence $soc_{\X} (C) \subset soc_{\X} (\oplus_{i\in I} \bar{M_i} ) $, and the intersection $C\cap Im\ g=0$ implies $C=0$. Therefore $g$ is surjective, hence an automorphism. Thus $f$ is left minimal, and therefore an $\X^{\perp}$-envelope of $M$. 
\\%
(b) If $M$ is $\Ext$-injective in $\Filt{\X}$, then $M= M_{\X}(\infty)$ as constructed in Definition \ref{construction_universal_sequence}. The claim now follows from part (a) together with Theorem \ref{envelope}.
\end{proof}
\begin{corollary} \label{countable_X-socle_filtration}
	Let $\X$ be a semibrick with $\pd{\X}\leq 1$ and without $\X$-divisible objects. Then the $\X$-socle filtration of every module $M$ in $\Filt{\X}$ has finite or countable length.
\end{corollary}
\begin{proof}
	Every module in $\Filt{\X}$ admits an $\X^{\perp}$-envelope, which by the theorem above  is isomorphic to a direct sum of Pr\"ufer modules. Each Pr\"ufer module has, by Proposition \ref{quotient_universal_sequences}(b), a countable $\X$-socle filtration. Hence, by  Proposition \ref{closed_under_cokernel}(b), every submodule in $\Filt{\X}$ also has a finite or countable $\X$-socle filtration.
\end{proof}
	If $R$ is a tame hereditary algebra, then all Pr\"ufer modules are algebraically compact. An elegant proof of this fact (see \cite[p. 46]{10.1007/978-3-0348-8426-6_1}) uses that a Pr\"ufer module, considered as a module over its endomorphism ring, is  artinian. 
	It is, however, not known whether Pr\"ufer modules as defined in this paper are algebraically compact.	
%
%
%

%
\end{document}